\documentclass[12pt]{conm-p-l}
\usepackage{epic}
\usepackage{eepic}
\usepackage{amsfonts} 
\usepackage{amssymb}
\usepackage[all]{xy}
\usepackage{bbm}
\usepackage{graphicx}
\usepackage{url}
\usepackage{color}

\newlength{\rememberparindent}
\newtheorem{theorem1}{Theorem}
\newtheorem{corollary1}[theorem1]{Corollary}
\newtheorem{proposition1}[theorem1]{Proposition}
\newtheorem{lemma1}[theorem1]{Lemma}
\newtheorem{problem1}{Problem}[section]
\newtheorem{conjecture1}[problem1]{Conjecture}
\newtheorem{question1}[problem1]{Question}

\theoremstyle{definition}
\newtheorem{definition1}[theorem1]{Definition}
\newtheorem{example1}[theorem1]{Example}

\theoremstyle{remark}
\newtheorem{remark1}[theorem1]{Remark}

\newenvironment{theorem}{\setlength{\parindent}{0pt}\begin{theorem1}}{\end{theorem1}\setlength{\parindent}{\rememberparindent}}

\newenvironment{lemma}{\setlength{\parindent}{0pt}\begin{lemma1}}{\end{lemma1}\setlength{\parindent}{\rememberparindent}}

\newcommand{\Z}{{\mathbb Z}}
\newcommand{\Q}{{\mathbb Q}}
\newcommand{\C}{{\mathbb C}}
\newcommand{\one}{{\mathbbm 1}}
\newcommand{\nc}{{\mathcal N}}

\newcommand{\ind}{\operatorname{ind}}
\newcommand{\suchthat}{\ : \ }

\renewcommand{\epsilon}{\varepsilon}
\newcommand{\gf}{G}
\newcommand{\T}{{\mathcal T}}

\begin{document}
\title[Polar decomposition and Brion's theorem.]{Polar decomposition
  and Brion's theorem.}
\author{Christian Haase}
\email{haase@math.duke.edu}
\address{Department of Mathematics \\
  Duke University \\
  Durham, NC 27708--0320 \\
  USA}
\thanks{Work supported by NSF grant DMS--0200740}

\begin{abstract}
  In this note we point out the relation between Brion's formula for
  the lattice point generating function of a convex polytope in terms
  of the vertex cones~\cite{BrionTangent} on the one hand, and the
  polar decomposition \`a la
  Lawrence/Varchenko~\cite{lawrenceVolume,varchenko} on the other.
  We then go on to prove a version of polar decomposition for
  non-simple polytopes.
\end{abstract}

\maketitle
\thispagestyle{empty}
\section{Introduction}\noindent
The stars of this note are two formulas that express convex polytopes
in terms of cones. On the one hand, Brion's theorem
(equation~\eqref{eq:BrionGF} below) expresses lattice
points in a polytope in terms of tangent cones at
vertices~\cite{BrionTangent}. On the other hand, given a simple
polytope and a generic objective function, one can write the indicator
function of a simple polytope as a combination of polarized
tangent cones at the vertices (equation~\eqref{eq:LawrenceIndicator}).
We refer to this formula as polar
decomposition~\cite{lawrenceVolume,varchenko}.

After we review these formulas in this introduction, we point out in
section~\ref{sec:LthusB} some simple relations between the two
formulas which, to our knowledge, have not appeared in the literature
before. In particular, Brion's theorem (for simple polytopes) follows
from polar decomposition.
The last section is devoted to the formulation and proof of polar
decomposition for non-simple polytopes\footnote{Put together,
  sections~\ref{sec:LthusB} and \ref{sec:non-simple-lawrence} yield
  another proof of Brion's theorem.}.

This is by no means an attempt to survey what is scattered throughout
the literature\footnote{See~\cite[p.~346]{BarvinokBook} and the
  references therein for Brion's theorem,
  \cite{lawrenceVolume,varchenko} for polar decomposition, and
  \cite{Brianchon,Gram,Shephard} for the Brianchon-Gram formula below.}
about decompositions of polytopes into cones. In fact, Matthias Beck
made me aware of such a survey in the making~\cite{conicSurvey}. We
will concentrate on what we think is new. We will not state the most
general form of the results, but leave generalizations (e.g., weighted
versions) as exercises for the ambitious reader.
\subsection*{{\normalfont \em Acknowledgments}}
This paper has its origins in Jos\'e Agapito's inspiring talk about
his weigh\-ted version of polar decomposition for simple
polytopes~\cite{pepe}. After the talk, I asked the question whether
one can obtain an analogous formula for non-simple polytopes by
considering simple deformations. I got intrigued by this question,
and, after computing several examples, was convinced that there is
such a theorem out there.

I was lucky to get hold of a preliminary draft of the above mentioned
survey~\cite{conicSurvey}. Not only did it take away the pressure to
try and write a comprehensive article about conic decompositions of
polytopes, but it also made me aware of references and view points
that I did not know about before. Much of the presentation in the
present paper is influenced by (if not stolen from) this draft.
Finally, I want to thank the referee who helped to clarify the
exposition.
\subsection{Basic definitions and the Brianchon-Gram formula}
We will consider rational convex polytopes $P$ and
polyhedral cones $C$ in $\Q^d$, where $\Q^d$ is endowed with the fixed
lattice $\Z^d$. We will assume that polytopes and cones are
full-dimensional.
For standard polytope definitions and notation we refer
to~\cite{Ziegler}.
To a subset $S \subseteq \Q^d$ we assign two objects:
\begin{list}{-}{
    \setlength{\topsep}{1mm}
    \setlength{\parsep}{\topsep}
    \setlength{\leftmargin}{\parindent}
    \setlength{\labelwidth}{1em}}
\item[-] The indicator function $\one_S \colon \Q^d \rightarrow \Z$
  vanishes outside $S$, and is one along $S$.
\item[-] The generating function (of the lattice points) is the formal
  Laurent power series $\gf_S = \sum_{m \in S \cap \Z^d} z^m \in
  \C[[z_1^{\pm 1}, \ldots, z_d^{\pm   1}]]$, where we write $z^m =
  z_1^{m_1} \cdot \ldots \cdot z_d^{m_d}$.
If $S=P$ is a polytope, then $\gf_P$ is a Laurent polynomial; when
$S=C$ is a cone, then $\gf_C$ is a rational function\footnote{We are
  glossing over significant technicalities here and later on when we
  pass from formal power series to rational functions. For a first
  reading, we prefer to skip this step, and refer
  to~\cite{BarvinokBook,bv} for a careful treatment. The experts
  should be familiar with the standard arguments, anyway.}.
\end{list}
The
essential features can be illustrated on the line. For $a < b \in \Z$,
we have $\gf_{[a,b]}(x) = x^a+\ldots+x^b = \frac{x^a-x^{b+1}}{1-x}$,
and $\gf_{[a,\infty)}(x) = x^a+\ldots = \frac{x^a}{1-x}$.

The mother of all conic decomposition theorems is the Brianchon-Gram
formula. According to Shephard~\cite{Shephard},
Brianchon~\cite{Brianchon} and Gram~\cite{Gram} independently proved
the $d=3$  case in 1837(!) and 1874 respectively.
In 1927 Sommerville~\cite{Sommerville} published a proof for
general $d$, which was corrected in the 1960's by
Gr\"unbaum~\cite[\S14.1]{Gruenbaum}.

\begin{figure}[hbtp]
  \centering
\begin{picture}(0,0)%
\includegraphics{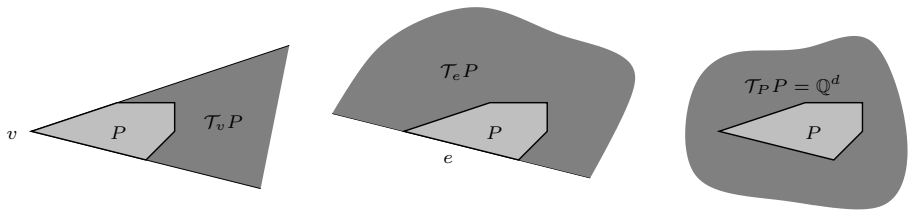}%
\end{picture}%
\setlength{\unitlength}{1184sp}%
\begingroup\makeatletter\ifx\SetFigFont\undefined%
\gdef\SetFigFont#1#2#3#4#5{%
  \reset@font\fontsize{#1}{#2pt}%
  \fontfamily{#3}\fontseries{#4}\fontshape{#5}%
  \selectfont}%
\fi\endgroup%
\begin{picture}(18870,4239)(1876,-4998)
\put(6001,-3286){\makebox(0,0)[lb]{\smash{\SetFigFont{7}{8.4}{\rmdefault}{\mddefault}{\updefault}{\color[rgb]{0,0,0}$\T_vP$}%
}}}
\put(1876,-3511){\makebox(0,0)[lb]{\smash{\SetFigFont{7}{8.4}{\rmdefault}{\mddefault}{\updefault}{\color[rgb]{0,0,0}$v$}%
}}}
\put(4051,-3511){\makebox(0,0)[lb]{\smash{\SetFigFont{7}{8.4}{\rmdefault}{\mddefault}{\updefault}{\color[rgb]{0,0,0}$P$}%
}}}
\put(18601,-3511){\makebox(0,0)[lb]{\smash{\SetFigFont{7}{8.4}{\rmdefault}{\mddefault}{\updefault}{\color[rgb]{0,0,0}$P$}%
}}}
\put(11926,-3511){\makebox(0,0)[lb]{\smash{\SetFigFont{7}{8.4}{\rmdefault}{\mddefault}{\updefault}{\color[rgb]{0,0,0}$P$}%
}}}
\put(11026,-4036){\makebox(0,0)[lb]{\smash{\SetFigFont{7}{8.4}{\rmdefault}{\mddefault}{\updefault}{\color[rgb]{0,0,0}$e$}%
}}}
\put(10951,-2236){\makebox(0,0)[lb]{\smash{\SetFigFont{7}{8.4}{\rmdefault}{\mddefault}{\updefault}{\color[rgb]{0,0,0}$\T_eP$}%
}}}
\put(17326,-2536){\makebox(0,0)[lb]{\smash{\SetFigFont{7}{8.4}{\rmdefault}{\mddefault}{\updefault}{\color[rgb]{0,0,0}$\T_PP=\Q^d$}%
}}}
\end{picture}
  \caption{
    Tangent cones 
    at vertex $v$, edge $e$, and at $P$.}
  \label{fig:tangentCone}
\end{figure}

Let $P \subset \Q^d$ be a convex polytope. For a face $F$ of $P$, we
define the tangent cone (of $P$ at $F$) by
\begin{equation*}
  \T_FP = \{ f+x \in \Q^d \suchthat f \in F \text{ and } f + \epsilon
  x \in P \text{ for some } \epsilon > 0 \}
\end{equation*}
(See Figure~\ref{fig:tangentCone}.)
Then $\one_P$ is the alternating sum of all tangent cones:
\begin{equation} \label{eq:Gram}
  \one_P = \sum_{F \preceq P} (-1)^{\dim F} \ \one_{\T_FP} .
\end{equation}
At every point $x \in P$, the right hand side computes the Euler
characteristic of $P$ (and $P$ is contractible). While outside $P$, we
have to subtract the Euler characteristic of the subcomplex that is
visible from $x$ (which again, is contractible, actually
shellable). (Compare~\cite{Shephard}.)
%
\subsection{Brion's formula}
Brion~\cite{BrionTangent} proved his decomposition using a
Riemann-Roch type formula. In the mean time, more elementary proofs
have been
given~\cite{BarvinokExponential,bv,LawrenceValuations,PukhlikovKhovanskii}.
The easiest way to formulate Brion's formula is in terms of
generating functions. We can express the generating function for $P$
as the sum of the generating functions of the tangent cones at
vertices.
\begin{equation} \label{eq:BrionGF}
  \gf_P = \sum_{v \text{ vertex of } P} \gf_{\T_vP}
\end{equation}
In the example, $\T_a[a,b] = [a,\infty)$, and $\T_b[a,b] = (-\infty,b]$
with generating functions $\gf_{[a,\infty)} = \frac{x^a}{1-x}$ and
$\gf_{(-\infty,b]} = \frac{(1/x)^{-b}}{1-1/x} =
-\frac{x^{b+1}}{1-x}$. As predicted, they sum to $\gf_{[a,b]}$.

This version of Brion's formula is implied by the following two
results. (Compare~\cite{BarvinokBook}.)
\begin{theorem} \label{thm:BrionIndicator}
  The function $\one_P - \sum_{v \text{ vertex of } P} \one_{\T_vP}$
  is a linear combination of indicator functions of cones that contain
  affine lines.
\end{theorem}
This follows immediately from~\eqref{eq:Gram}.
It is also implied by the considerations in Sections~\ref{sec:LthusB}
and~\ref{sec:non-simple-lawrence}.
In our example, $\one_{[a,b]} = \one_{[a,\infty)} + \one_{(\infty,b]}
- \one_\Q$.
\begin{lemma}
  Let $C \subseteq \Q^d$ be a cone that contains an affine line. Then
  $\gf_C=0$.
\end{lemma}
The idea of the proof is to decompose the monoid $C \cap \Z^d$ into
a direct sum $L \oplus R$ of a one dimensional lattice $L$ and a
complement $R$. Then $\gf_C$ can be written as $\gf_R \gf_L$, and
$\gf_L=0$\footnote{Glossing over alert: look under the
  rug~\cite{BarvinokBook}!}.
\subsection{(Simple) polar decomposition}
By flipping the edge vectors emanating from each vertex of $P$ in a
systematic way, we can get the so-called polar decomposition, which
expresses the characteristic function of a convex simple polytope in
terms of the characteristic functions of the cones supported at the
vertices of P (no cones with straight lines needed this time). This
was first obtained by Varchenko~\cite{varchenko} and
Lawrence~\cite{lawrenceVolume} independently. Lawrence used
Brianchon-Gram's theorem together with the principle of
inclusion-exclusion to derive the polar decomposition theorem. Then
Karshon, Sternberg and Weitsman obtained a weighted version of this
decomposition by assigning, in a consistent way, particular weights to
the lattice points in the polytope and in the
cones~\cite{KarshonSternbergWeitsman}. Their work was also motivated
by methods in differential and algebraic geometry. Finally, using the
same source of motivation, Agapito~\cite{pepe} gave a more general
weighted version of the polar decomposition theorem for simple
polytopes that includes Lawrence/Varchenko and
Karshon-Sternberg-Weitsman versions as particular cases.

Polar decomposition is like \allowbreak Morse theory for simple
polytopes. We sweep over the polytope and build it up from local
contributions, critical point (vertex) by critical point.
As input we use a generic\footnote{Generic means here that $\xi$ is
  not constant on any edge of $P$. In
  Section~\ref{sec:non-simple-lawrence} we ask a little bit more.}
linear functional $\xi \in (\Q^d)^*$. Using $\xi$, we define the
polarized tangent cones as follows.

Because $P$ is simple, at every vertex $v$, the tangent cone $\T_vP$
is generated by $d$ linearly independent directions $t_1, \ldots, t_d$
(remember $P$ is full-dimensional).
\begin{equation*}
  \T_vP = v + \sum_{i=1}^d \Q_{\ge 0} t_i
\end{equation*}
Now the polarized tangent cone (with respect to $\xi$) is
\begin{equation*}
  \T^\xi_vP = v + \sum_{\xi(t_i) > 0} \Q_{\ge 0} t_i + \sum_{\xi(t_i)
    < 0} \Q_{< 0} t_i .
\end{equation*}
\begin{figure}[htbp]
  \begin{center}
\begin{picture}(0,0)%
\includegraphics{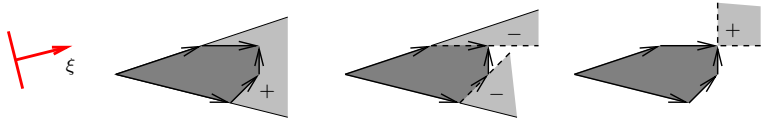}%
\end{picture}%
\setlength{\unitlength}{1184sp}%
\begingroup\makeatletter\ifx\SetFigFont\undefined%
\gdef\SetFigFont#1#2#3#4#5{%
  \reset@font\fontsize{#1}{#2pt}%
  \fontfamily{#3}\fontseries{#4}\fontshape{#5}%
  \selectfont}%
\fi\endgroup%
\begin{picture}(15837,2445)(97,-9673)
\put(1351,-8686){\makebox(0,0)[lb]{\smash{\SetFigFont{7}{8.4}{\rmdefault}{\mddefault}{\updefault}{\color[rgb]{0,0,0}$\xi$}%
}}}
\put(5401,-9211){\makebox(0,0)[lb]{\smash{\SetFigFont{7}{8.4}{\rmdefault}{\mddefault}{\updefault}{\color[rgb]{0,0,0}$+$}%
}}}
\put(10576,-8011){\makebox(0,0)[lb]{\smash{\SetFigFont{7}{8.4}{\rmdefault}{\mddefault}{\updefault}{\color[rgb]{0,0,0}$-$}%
}}}
\put(15151,-7936){\makebox(0,0)[lb]{\smash{\SetFigFont{7}{8.4}{\rmdefault}{\mddefault}{\updefault}{\color[rgb]{0,0,0}$+$}%
}}}
\put(10201,-9286){\makebox(0,0)[lb]{\smash{\SetFigFont{7}{8.4}{\rmdefault}{\mddefault}{\updefault}{\color[rgb]{0,0,0}$-$}%
}}}
\end{picture}
    \caption{Polarized tangent cones.}
    \label{fig:tangentConePlus}
  \end{center}
\end{figure}
This is a locally closed cone all whose points are $\xi$-higher than
$v$. (Compare Figure~\ref{fig:tangentConePlus}.) Dually, if $n_1,
\ldots, n_d \in (\Q^d)^*$ are the inner facet normals to $P$ at $v$,
and $\xi=\sum \alpha_i n_i$, then $\T^\xi_vP$ can be defined by the
inequalities
\begin{equation} \label{eq:TxiPineq}
  n_i(x) \ge n_i(v) \text{ for } \alpha_i > 0, \text{ and } n_i(x) <
  n_i(v) \text{ for } \alpha_i < 0 .
\end{equation}
Define the index $\ind_\xi(v)$ of the vertex $v$ as the number of
$\xi$-negative edge directions $t_i$, or equivalently, as the number
of negative coefficients $\alpha_i$. Then we can write 
$P$ as the signed sum of 
polarized tangent cones.
\begin{equation} \label{eq:LawrenceIndicator}
  \one_P = \sum_{v \text{ vertex of } P} (-1)^{\ind_\xi(v)}
  \one_{\T^\xi_vP}
\end{equation}
Again, we can formulate a generating function version of this result:
just replace all $\one$'s by $\gf$'s\footnote{More glossing over.}.

Equation~\eqref{eq:LawrenceIndicator} follows from \eqref{eq:Gram} if
we group all faces according to where they achieve their
$\xi$-maximum. Then it remains to check that
\begin{equation} \label{eq:rearrange}
  (-1)^{\ind_\xi(v)} \one_{\T^\xi_vP} = \sum_{
    \substack{v \preceq F \preceq P \\[.5mm] \xi(F) \le \xi(v)}
  }
  (-1)^{\dim F} \ \one_{\T_FP} .
\end{equation}
There is a weighted version of
equation~\eqref{eq:LawrenceIndicator}~\cite{pepe}. The weighted
indicator function of a polytope $P$ is the function $\one^w_P \colon
\Q^d \rightarrow \Z[z]$ which vanishes outside $P$, and takes the value
$z^k$ along the relative interior of codimension $k$
faces.\footnote{The formula in~\cite{pepe} is stated using the
  substitution $z=\frac{1}{1+y}$.} One could even introduce a
different variable for every facet of $P$, and assign their product to
the corresponding intersection.
Now the same formula holds if we modify the indicator functions of the
polarized tangent cones as follows: 
A face of (the closure of) $\T^\xi_vP$ is defined as the set of points
that satisfy equality for some of the
inequalities~\eqref{eq:TxiPineq}. The value of $\one^w_{\T^\xi_vP}$
along a face defined by $k_+$ equalities $n_i(x) = n_i(v)$ for
$\alpha_i > 0$ and $k_-$ equalities for $\alpha_i < 0$ is the
polynomial $z^{k_+}(1-z)^{k_-} \in \Z[z]$. One recovers the unweighted
version for $z=1$, and for $z=0$ one obtains a decomposition of the
interior of $P$.
\section{This follows from that\label{sec:LthusB}}
\noindent
In this section we show the relation between the two formulas.
Essentially, for the indicator functions, Brion's formula (in the
simple case) is polar decomposition modulo cones that contain lines.
%
\begin{lemma}
  If $P$ is a simple polytope and $v$ is a vertex, then the different
  polarized tangent cones for various $\xi$ partition $\Q^d$. Moreover,
  for $\xi_1, \xi_2 \in (\Q^d)^*$, $(-1)^{\ind_{\xi_1}(v)}
  \one_{\T^{\xi_1}_vP}$ and $(-1)^{\ind_{\xi_2}(v)}
  \one_{\T^{\xi_2}_vP}$ are equivalent modulo cones that contain
  lines.
\end{lemma}
\begin{proof}
  If $n_i(x) \ge n_i(v)$ for $i = 1, \ldots, d$ are the inequalities
  of the facets incident to $v$, the hyperplanes $n_i(x)=n_i(v)$
  subdivide $\Q^d$ into orthant cones.
  Every point $x \in \Q^d$ belongs to exactly one polarized tangent
  cone according to the signs of the $n_i(x)$.

  We can get from every polarized tangent cone to every other by
  successively flipping inequalities. Now the sum of two adjacent
  polarized tangent cones is defined by $d-1$ (strict and non-strict)
  inequalities. Thus it is a cone that contains a line.
\end{proof}
For a simple polytope, the indicator function version of Brion's
formula (Theorem~\ref{thm:BrionIndicator}) 
can be derived from
equation~\eqref{eq:LawrenceIndicator} (modulo cones that contain
lines).
As Theorem~\ref{thm:BrionIndicator} does not specify these cones,
there is no converse. On the other hand, generating functions do not
see cones that contain lines. So for generating functions, the two
formulas actually coincide.

The same considerations carry through for the weighted indicator
functions. One thus obtains a weighted Brion's formula. We leave this
as an exercise for the reader.
\section{Non-simple polar decomposition\label{sec:non-simple-lawrence}}
\noindent
We want to generalize polar decomposition to non-simple polytopes. We
will compute the local contribution from a simple deformation of the
vertex in question, and prove that the result does not depend on the
chosen deformation.

As there are not too many non-simple polytopes in dimension $\le 2$,
the new running example will be the pyramid with the five vertices
$(0,0,0)$, and $(\pm 1, \pm 1, 1)$.
\begin{figure}[htbp]
  \centering
\setlength{\unitlength}{0.00025000in}
\begingroup\makeatletter\ifx\SetFigFont\undefined%
\gdef\SetFigFont#1#2#3#4#5{%
  \reset@font\fontsize{#1}{#2pt}%
  \fontfamily{#3}\fontseries{#4}\fontshape{#5}%
  \selectfont}%
\fi\endgroup%
{\renewcommand{\dashlinestretch}{30}
\begin{picture}(3624,3039)(0,-10)
\path(12,1812)(2412,1512)(3612,2712)
        (1212,3012)(12,1812)
\path(12,1812)(1812,12)(2412,1512)
\path(1812,12)(3612,2712)
\dashline{120.000}(1812,12)(1212,3012)
\end{picture}
}
  \caption{The simplest non-simple polytope.}
  \label{fig:3deg}
\end{figure}
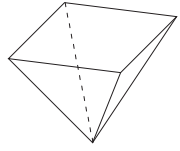

Let $n_1, \ldots, n_N$ be the inner facet normals to $P$ at $v$.
They generate the rays of the normal cone $\nc_v$. A virtual
deformation of the vertex $v$ is a regular
triangulation\footnote{Compare~\cite[\S~14.3]{LeeHandBook}. Moving the
  facets according to the values of the convex function at the $n_i$
  would yield an actual simple deformation of the vertex $v$.}
$\Delta$ of $\nc_v$. That is a face-to-face subdivision of $\nc_v$
into simplicial cones $\sigma_1, \ldots, \sigma_M$, so that there is a
convex piecewise linear function on $\nc_v$ with domains of linearity
the $\sigma_i$.

\begin{figure}[htbp]
  \centering
  \raisebox{7mm}{$(\Q^d)^*$ :} \qquad
\setlength{\unitlength}{0.00016667in}
\begingroup\makeatletter\ifx\SetFigFont\undefined%
\gdef\SetFigFont#1#2#3#4#5{%
  \reset@font\fontsize{#1}{#2pt}%
  \fontfamily{#3}\fontseries{#4}\fontshape{#5}%
  \selectfont}%
\fi\endgroup%
{\renewcommand{\dashlinestretch}{30}
\begin{picture}(4374,4287)(0,-10)
\path(387,2937)(3987,2487)
\path(2187,12)(987,2412)
\path(1188.246,2177.213)(987.000,2412.000)(1054.082,2110.131)
\path(2187,12)(4362,3012)
\path(4246.631,2725.094)(4362.000,3012.000)(4125.189,2813.139)
\path(2187,12)(12,3537)
\path(233.360,3321.072)(12.000,3537.000)(105.704,3242.306)
\dashline{120.000}(2187,12)(3087,3612)
\path(3087,3612)(3237,4212)
\path(3237.000,3902.767)(3237.000,4212.000)(3091.479,3939.147)
\path(1287,1812)(3987,2487)(3087,3612)
        (387,2937)(1287,1812)
\put(4062,1887){\makebox(0,0)[lb]{\smash{{{\SetFigFont{10}{12.0}{\rmdefault}{\mddefault}{\updefault}$n_3$}}}}}
\put(312,3462){\makebox(0,0)[lb]{\smash{{{\SetFigFont{10}{12.0}{\rmdefault}{\mddefault}{\updefault}$n_4$}}}}}
\put(2112,3912){\makebox(0,0)[lb]{\smash{{{\SetFigFont{10}{12.0}{\rmdefault}{\mddefault}{\updefault}$n_2$}}}}}
\put(162,1662){\makebox(0,0)[lb]{\smash{{{\SetFigFont{10}{12.0}{\rmdefault}{\mddefault}{\updefault}$n_1$}}}}}
\end{picture}
}
 \qquad
\setlength{\unitlength}{0.00016667in}
\begingroup\makeatletter\ifx\SetFigFont\undefined%
\gdef\SetFigFont#1#2#3#4#5{%
  \reset@font\fontsize{#1}{#2pt}%
  \fontfamily{#3}\fontseries{#4}\fontshape{#5}%
  \selectfont}%
\fi\endgroup%
{\renewcommand{\dashlinestretch}{30}
\begin{picture}(4374,4287)(0,-10)
\path(2187,12)(987,2412)
\path(1188.246,2177.213)(987.000,2412.000)(1054.082,2110.131)
\path(2187,12)(4362,3012)
\path(4246.631,2725.094)(4362.000,3012.000)(4125.189,2813.139)
\path(2187,12)(12,3537)
\path(233.360,3321.072)(12.000,3537.000)(105.704,3242.306)
\dashline{120.000}(2187,12)(3087,3612)
\path(3087,3612)(3237,4212)
\path(3237.000,3902.767)(3237.000,4212.000)(3091.479,3939.147)
\path(1287,1812)(3987,2487)(3087,3612)
        (387,2937)(1287,1812)
\path(1287,1812)(3087,3612)
\put(4062,1887){\makebox(0,0)[lb]{\smash{{{\SetFigFont{10}{12.0}{\rmdefault}{\mddefault}{\updefault}$n_3$}}}}}
\put(312,3462){\makebox(0,0)[lb]{\smash{{{\SetFigFont{10}{12.0}{\rmdefault}{\mddefault}{\updefault}$n_4$}}}}}
\put(2112,3912){\makebox(0,0)[lb]{\smash{{{\SetFigFont{10}{12.0}{\rmdefault}{\mddefault}{\updefault}$n_2$}}}}}
\put(162,1662){\makebox(0,0)[lb]{\smash{{{\SetFigFont{10}{12.0}{\rmdefault}{\mddefault}{\updefault}$n_1$}}}}}
\end{picture}
}
 \\[3mm]
\raisebox{7mm}{$\Q^d$ :} \qquad
\setlength{\unitlength}{0.00016667in}
\begingroup\makeatletter\ifx\SetFigFont\undefined%
\gdef\SetFigFont#1#2#3#4#5{%
  \reset@font\fontsize{#1}{#2pt}%
  \fontfamily{#3}\fontseries{#4}\fontshape{#5}%
  \selectfont}%
\fi\endgroup%
{\renewcommand{\dashlinestretch}{30}
\begin{picture}(4224,3639)(0,-10)
\path(12,1812)(2412,1512)(4212,3312)
        (1812,3612)(12,1812)
\path(12,1812)(1812,12)(2412,1512)
\dashline{120.000}(2412,612)(1812,3612)
\path(2412,612)(4212,3312)
\path(1812,12)(2412,612)
\end{picture}
}
 \qquad
\setlength{\unitlength}{0.00016667in}
\begingroup\makeatletter\ifx\SetFigFont\undefined%
\gdef\SetFigFont#1#2#3#4#5{%
  \reset@font\fontsize{#1}{#2pt}%
  \fontfamily{#3}\fontseries{#4}\fontshape{#5}%
  \selectfont}%
\fi\endgroup%
{\renewcommand{\dashlinestretch}{30}
\begin{picture}(4224,3114)(0,-10)
\path(12,1887)(3012,1512)(4212,2712)
        (1212,3087)(12,1887)
\path(2412,12)(4212,2712)
\dashline{120.000}(1812,87)(1212,3087)
\path(12,1887)(1812,87)(2412,12)(3012,1512)
\end{picture}
}
  \caption{Triangulations and the corresponding deformations.}
  \label{fig:deform}
\end{figure}
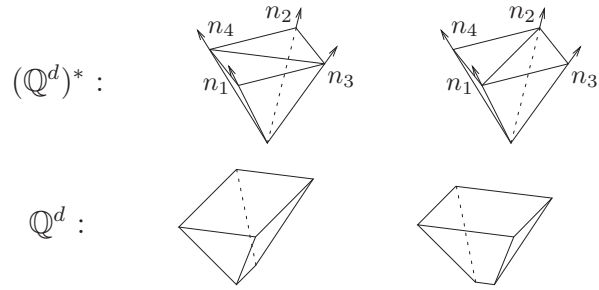
In the example, $v = (0,0,0)$ is the only non-simple vertex. The
normal cone is generated by $n_1=(1,0,1)$, $n_2=(-1,0,1)$,
$n_3=(0,1,1)$, and $n_4=(0,-1,1)$. It has two triangulations, both of
which are regular:
$\Delta_1$ with maximal simplices spanned by $n_1,n_3,n_4$ and
$n_2,n_3,n_4$, and $\Delta_2$ with maximal simplices spanned by
$n_1,n_2,n_3$ and $n_1,n_2,n_4$. The corresponding deformations are
sketched in Figure~\ref{fig:deform}.

Now the tangent cone $\T_vP$ can be written as the intersection of
simple cones defined by the $\sigma_i$:
\begin{equation} \label{eq:nonSimpleTangent}
  \T_vP = \bigcap_i \T_{\sigma_i} \text{ for } \T_{\sigma_i} = \{ x
  \in \Q^d \,:\, n_j(x) \ge n_j(v) \text{ for all } n_j \in
  \sigma_i \} .
\end{equation}
For generic\footnote{Now, generic means that $\xi$ is not constant
  on any ray of any of the $\T_{\sigma_i}$'s, i.e., $\xi$ does not lie
  on any hyperplane used in the triangulation.} $\xi \in (\Q^d)^*$, we
write $\xi$ in terms of the $n_j \in \sigma_i$, and flip those
inequalities in~\eqref{eq:nonSimpleTangent} where $\xi$ has negative
coefficients. Then we compute the index $\ind_\xi(\sigma_i)$ as the
number of flipped inequalities, and define the local contribution at
$v$ as:
\begin{equation*} 
  \one^\xi_{\Delta,v} = \sum_i (-1)^{\ind_\xi(\sigma_i)}
  \one_{\T^\xi_v \T_{\sigma_i}}
\end{equation*}
The cool thing is that $\one^\xi_{\Delta,v} = \one^\xi_v$ does not
depend on the triangulation. So it really is a {\em local\/}
contribution.
\begin{theorem}
  Let $P \subset \Q^d$ be any polytope, and let $\xi \in (\Q^d)^*$ be
  generic. Then $\one^\xi_{\Delta,v} = \one^\xi_v$ does not depend on
  the choice of a regular triangulation $\Delta$ of the inner normal
  cone $\nc_v$ at vertex $v$.
  Moreover, 
  \begin{equation} \label{eq:decomposition}
    \one_P(x) = \sum_{v \text{ vertex of } P} \one^\xi_v(x).
  \end{equation}
\end{theorem}
The $\Delta$-invariance can be shown, e.g., using the fact that all
regular triangulations are connected by
flips~\cite[\S~14.6]{LeeHandBook}. In this note, however, we follow a
different strategy. First, we will show the decomposition formula
\eqref{eq:decomposition} for compatible choices of triangulations
(Lemma \ref{lemma:compatible}). Then, we observe that under relatively
weak conditions, all decompositions are the same
(Lemma~\ref{lemma:local}).

\begin{figure}[htbp]
  \centering
\begin{picture}(0,0)%
\includegraphics{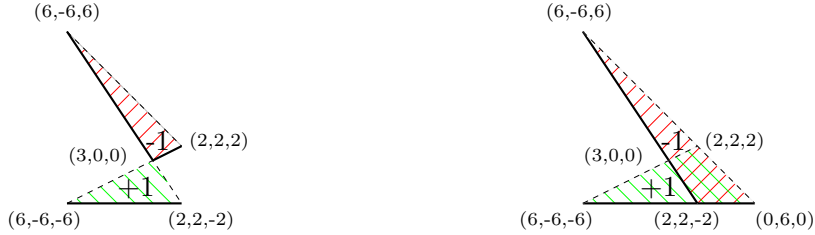}%
\end{picture}%
\setlength{\unitlength}{1184sp}%
\begingroup\makeatletter\ifx\SetFigFont\undefined%
\gdef\SetFigFont#1#2#3#4#5{%
  \reset@font\fontsize{#1}{#2pt}%
  \fontfamily{#3}\fontseries{#4}\fontshape{#5}%
  \selectfont}%
\fi\endgroup%
\begin{picture}(14488,4713)(1157,-5710)
\put(2251,-5011){\makebox(0,0)[lb]{\smash{\SetFigFont{11}{13.2}{\rmdefault}{\mddefault}{\updefault}{\color[rgb]{0,0,0}+1}%
}}}
\put(3301,-5611){\makebox(0,0)[lb]{\smash{\SetFigFont{7}{8.4}{\rmdefault}{\mddefault}{\updefault}{\color[rgb]{0,0,0}(2,2,-2)}%
}}}
\put(3751,-3961){\makebox(0,0)[lb]{\smash{\SetFigFont{7}{8.4}{\rmdefault}{\mddefault}{\updefault}{\color[rgb]{0,0,0}(2,2,2)}%
}}}
\put(2851,-4036){\makebox(0,0)[lb]{\smash{\SetFigFont{11}{13.2}{\rmdefault}{\mddefault}{\updefault}{\color[rgb]{0,0,0}-1}%
}}}
\put(2476,-4261){\makebox(0,0)[rb]{\smash{\SetFigFont{7}{8.4}{\rmdefault}{\mddefault}{\updefault}{\color[rgb]{0,0,0}(3,0,0)}%
}}}
\put(1501,-5611){\makebox(0,0)[rb]{\smash{\SetFigFont{7}{8.4}{\rmdefault}{\mddefault}{\updefault}{\color[rgb]{0,0,0}(6,-6,-6)}%
}}}
\put(1201,-1261){\makebox(0,0)[b]{\smash{\SetFigFont{7}{8.4}{\rmdefault}{\mddefault}{\updefault}{\color[rgb]{0,0,0}(6,-6,6)}%
}}}
\put(15601,-5611){\makebox(0,0)[lb]{\smash{\SetFigFont{7}{8.4}{\rmdefault}{\mddefault}{\updefault}{\color[rgb]{0,0,0}(0,6,0)}%
}}}
\put(14551,-3961){\makebox(0,0)[lb]{\smash{\SetFigFont{7}{8.4}{\rmdefault}{\mddefault}{\updefault}{\color[rgb]{0,0,0}(2,2,2)}%
}}}
\put(13651,-4036){\makebox(0,0)[lb]{\smash{\SetFigFont{11}{13.2}{\rmdefault}{\mddefault}{\updefault}{\color[rgb]{0,0,0}-1}%
}}}
\put(13276,-4261){\makebox(0,0)[rb]{\smash{\SetFigFont{7}{8.4}{\rmdefault}{\mddefault}{\updefault}{\color[rgb]{0,0,0}(3,0,0)}%
}}}
\put(12226,-5611){\makebox(0,0)[rb]{\smash{\SetFigFont{7}{8.4}{\rmdefault}{\mddefault}{\updefault}{\color[rgb]{0,0,0}(6,-6,-6)}%
}}}
\put(14176,-5611){\makebox(0,0)[b]{\smash{\SetFigFont{7}{8.4}{\rmdefault}{\mddefault}{\updefault}{\color[rgb]{0,0,0}(2,2,-2)}%
}}}
\put(11926,-1261){\makebox(0,0)[b]{\smash{\SetFigFont{7}{8.4}{\rmdefault}{\mddefault}{\updefault}{\color[rgb]{0,0,0}(6,-6,6)}%
}}}
\put(13201,-5011){\makebox(0,0)[lb]{\smash{\SetFigFont{11}{13.2}{\rmdefault}{\mddefault}{\updefault}{\color[rgb]{0,0,0}+1}%
}}}
\end{picture}
  \caption{The local contribution from $\Delta_1$ and from $\Delta_2$.}
  \label{fig:flip}
\end{figure}

But let us first see how the
$\Delta$-invariance works out in our example.
If we choose $\xi=(4,2,0)$, 
then we compute for triangulation
$\Delta_1$ as follows: $\xi = 4 n_1 - n_3 - 3 n_4 = -4 n_2 + 3 n_3 +
n_4$. Thus the polarized cones are
\begin{align*}
  \T_v^\xi \T_{\sigma_1} &= \{ x \suchthat n_1(x) \ge 0, n_3(x) < 0,
  n_4(x) < 0 \} \text{ with index } 2 \text{ and } \\
  \T_v^\xi \T_{\sigma_2} &= \{ x \suchthat n_2(x) < 0, n_3(x) \ge 0,
  n_4(x) \ge 0 \} \text{ with index } 1.
\end{align*}
The trace of $\one_{\T_v^\xi \T_{\sigma_1}}-\one_{\T_v^\xi
  \T_{\sigma_2}}$ in the $\xi = 12$ plane is sketched in
Figure~\ref{fig:flip} on the left.
Similarly, for $\Delta_2$ we get
\begin{align*}
  \T_v^\xi \T_{\sigma_3} &= \{ x \suchthat n_1(x) \ge 0, n_2(x) < 0,
  n_3(x) \ge 0 \} \text{ with index } 1 \text{ and } \\
  \T_v^\xi \T_{\sigma_4} &= \{ x \suchthat n_1(x) \ge 0, n_2(x) < 0,
  n_4(x) < 0 \} \text{ with index } 2.
\end{align*}
The trace of $\one_{\T_v^\xi \T_{\sigma_4}}-\one_{\T_v^\xi
  \T_{\sigma_3}}$ in the $\xi = 12$ plane is sketched in
Figure~\ref{fig:flip} on the right. Observe that $\T_v^\xi
\T_{\sigma_3}$ and $\T_v^\xi \T_{\sigma_4}$ cancel out on the overlap
to yield the same contribution as $\T_v^\xi \T_{\sigma_1}$ and
$\T_v^\xi \T_{\sigma_2}$.
\newpage
\begin{lemma} \label{lemma:compatible}
  Let $P \subset \Q^d$ be any polytope, and let $\xi \in (\Q^d)^*$ be
  generic. Suppose $\Delta$ is a regular triangulation of the polar
  dual polytope $P^*$. It restricts to regular triangulations
  $\Delta_v$ of the normal cones $\nc_v$. Then
  \begin{equation} \label{eq:compatible}
    \one_P(x) = \sum_{v \text{ vertex of } P} \one^\xi_{\Delta_v,v}(x)
  \end{equation}
\end{lemma}
One can mimic the proof for the simple version
\eqref{eq:LawrenceIndicator} in order to show that the sum of the
$\one^\xi_{\Delta_v,v}$ over all vertices of $P$ does not depend on
$\xi$, and that for every $x \in \Q^d$ there is a suitable $\xi$ such
that equation~\eqref{eq:compatible}
is satisfied.

Now we actually use \eqref{eq:compatible} together with the following
lemma in order to show that $1^\xi_{\Delta,v} = 1^\xi_v$ is
independent of the triangulation.
\begin{lemma} \label{lemma:local}
  Let $P \subset \Q^d$ be any polytope, and let $\xi$ be a generic
  element of $(\Q^d)^*$. Suppose that there are two decompositions
  \begin{equation*}
    \one_P \ = \sum_{v \text{ vertex of } P} f_v \quad = \sum_{v
      \text{ vertex of } P} g_v
  \end{equation*}
  that are both
  {\em conic:\/}
  for every direction $t \in \Q^d \setminus 0$, the
  $f_v(v+\lambda t)$ and $g_v(v+\lambda t)$ are constant in $\lambda >
  0$, and
  {\em positive:\/} $f_v(v+t)=g_v(v+t)=0$ if $\xi(t) < 0$.
  Then $f_v = g_v$.
\end{lemma}
\begin{figure}[htbp]
  \begin{center}
\begin{picture}(0,0)%
\includegraphics{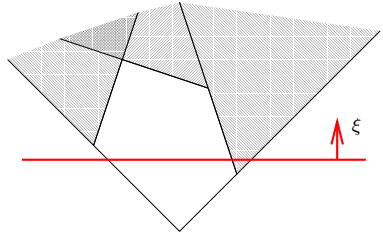}%
\end{picture}%
\setlength{\unitlength}{1184sp}%
\begingroup\makeatletter\ifx\SetFigFont\undefined%
\gdef\SetFigFont#1#2#3#4#5{%
  \reset@font\fontsize{#1}{#2pt}%
  \fontfamily{#3}\fontseries{#4}\fontshape{#5}%
  \selectfont}%
\fi\endgroup%
\begin{picture}(7824,4824)(3589,-8173)
\put(10801,-6061){\makebox(0,0)[lb]{\smash{\SetFigFont{7}{8.4}{\rmdefault}{\mddefault}{\updefault}{\color[rgb]{0,0,0}$\xi$}%
}}}
\end{picture}
  \end{center}
  \caption{Positive conic decomposition of $\one_P$ is unique.}
  \label{fig:posConic}
\end{figure}
\begin{proof}
  This proof illustrates what was meant by Morse theory earlier
  on. We sweep over the polytope bottom to top, and only at the
  vertices does something happen. That is, we compare the restrictions
  of the functions to the hyperplanes $\xi=c$.
  Order the vertices $v_1, \ldots, v_k$ in $\xi$-increasing order.
  For 
  $c \in [\xi(v_i),\xi(v_{i+1})[$,
  we have
  \begin{equation*}
    g_{v_i} \overset{(\mathbf{p})}{=} \sum_{j=i}^k
    g_{v_j} = \one_P - \sum_{j=1}^{i-1}
    g_{v_j} \overset{(\mathbf{c})}{=} \one_P - \sum_{j=1}^{i-1}
    f_{v_j} = \sum_{j=i}^k f_{v_j} \overset{(\mathbf{p})}{=} f_{v_i}
  \end{equation*}
  along $\xi=c$. For $(\mathbf{p})$ we use positivity of the $f_v$'s and
  $g_v$'s, and for $(\mathbf{c})$ we use induction, and the fact that
  the $f_v$'s and $g_v$'s are conic.
\end{proof}
Another incarnation of {\em the\/} positive conic decomposition is to
group the summands in the Brianchon-Gram formula~\eqref{eq:Gram} as we
did in equation~\eqref{eq:rearrange}, though it is less obvious that
this is positive.

\subsection{Homework. {\normalfont \em or\/} Why this is no good}
This is the wild speculation section. At the risk of exposing the full
extent of my ignorance, I ask a bunch of questions that I stumbled
over while compiling these notes.

If a good proof is one that makes us wiser~\cite{maninWiser}, then
this is a bad proof. It would be nicer to have a deformation
independent definition of the local contribution in terms of the
facet hyperplane arrangement at each vertex.
How much geometry is needed? Can we compute the local contribution in
purely combinatorial terms?

In some sense, the actual deformation of a vertex into several simple
vertices converges to the non-simple vertex, and the local contribution
is continuous. While this asks for a topology on the polytope algebra,
I believe that this is a more combinatorial question. There must be
a combinatorial deformation theory for (Eulerian?) posets.

Where are Gr\"obner bases? Often when regular triangulations show up,
the corresponding convex function determines a term order. Does the
result tell us something about commutative algebra?

\bibliographystyle{alpha}
\bibliography{alles}
\setlength{\parindent}{0pt}

\end{document}